\documentclass[12pt]{amsart}
\usepackage{amsmath,amsthm,amsfonts,amssymb}
\usepackage[all]{xy}

\DeclareMathOperator{\Aut}{Aut}
\DeclareMathOperator{\End}{End}
\DeclareMathOperator{\Fr}{Fr}

\DeclareMathOperator{\Hom}{Hom}

\DeclareMathOperator{\Mod}{Mod}

\DeclareMathOperator{\Id}{Id}

\newcommand{\C}{\mathbb C}
\newcommand{\CC}{\mathcal{C}}

\newcommand{\ot}{\otimes}

\newcommand{\B}{\mathcal{B}}
\newcommand{\PP}{\mathcal{P}}

\numberwithin{equation}{section}
\newtheorem{thm}[equation]{Theorem}
\newtheorem{theorem}[equation]{Theorem}

\newtheorem{corollary}[equation]{Corollary}
\newtheorem{lemma}[equation]{Lemma}

\theoremstyle{definition}

\newtheorem{rmk}[equation]{Remark}

\begin{document}

\title[Braid group representations]
{Braid group representations from twisted quantum doubles of finite groups}
\author{Pavel Etingof}
\address{Department of Mathematics,
   Massachusetts Institute of Technology,
  Cambridge, MA 02139, USA}
\email{etingof@math.mit.edu}
\author{Eric Rowell}
\address{Department of Mathematics,
    Texas A\&M University, College Station, TX 77843, USA}
   \email{rowell@math.tamu.edu}
\author{Sarah Witherspoon}
\email{sjw@math.tamu.edu}
\date{March 8, 2007}

\begin{abstract}
We investigate the braid group representations arising from categories
of representations of twisted quantum doubles of finite groups.
For these categories, we show that the resulting braid group
representations always factor through finite groups, in contrast
to the categories associated with quantum groups at roots of unity.
We also show that in the case of p-groups, the corresponding
pure braid group representations factor through a finite p-group,
which answers a question asked of the first author by V.\ Drinfeld.
\end{abstract}
\maketitle

\section{Introduction}

Any braided tensor category $\CC$ gives rise to finite dimensional representations of the braid group $\B_n$.  A natural problem is to determine the image of these representations.  This has been carried out to some extent for the braided tensor categories coming from quantum groups and polynomial link invariants at roots of unity \cite{FRW,FLW,GJ,J86,J89,LR,LRW}.  A basic question in this direction is: \emph{Is the image of the
representation of $\B_n$  a finite group?}  In the aforementioned papers the answer is typically ``no": Finite groups appear only in a few cases when the degree of the root of unity is small.

In this paper we consider the braid group representations associated to
the (braided tensor) categories $\Mod$-$D^{\omega}(G)$, where $D^{\omega}(G)$ is the
twisted quantum double of the finite group $G$.  We show (Theorem \ref{main}) that the
braid group images are \emph{always} finite.
We also answer in the affirmative (Theorem \ref{pgroup})
a question of Drinfeld: \emph{If $G$ is a $p$-group, is the image of the pure
  braid group $\PP_n$ also a $p$-group?}

The contents of the paper are as follows. In Section \ref{btc}
we record some definitions and basic results on braided
categories, and Section \ref{dg} is dedicated to the needed facts about $D^{\omega}(G)$.
Then we prove our main results in Section \ref{mt}.
The last section describes some open problems suggested by our work.

{\bf Acknowledgments.} P.E.\ is grateful to V.\ Drinfeld for a useful
discussion, and raising the question answered by Theorem \ref{pgroup}.
The work of P.E.\ was  partially supported by the NSF grant
 DMS-0504847. The work of S.W.\ was partially supported by the NSA
grant H98230-07-1-0038.

\section{Braided Categories and Braid Groups}\label{btc}
In this section we recall some facts about braided categories and derive
some basic consequences.  For more complete definitions the reader is referred to
either \cite{BK} or \cite{Tur1}.

The \textbf{braid group} $\B_n$ is defined
by generators $\beta_1,\ldots\beta_{n-1}$ satisfying the relations:
\begin{enumerate}
\item[(B1)] $\beta_i\beta_{i+1}\beta_i=\beta_{i+1}\beta_i\beta_{i+1}$ for $1\leq i\leq n-2,$
\item[(B2)] $\beta_i\beta_j=\beta_j\beta_i$ if $|i-j|\geq 2$.
\end{enumerate}

The kernel of the surjective homomorphism from $\B_n$ to the symmetric group $S_n$ given by $\beta_i\mapsto (i,i+1)$ is the \textbf{pure braid group} $\PP_n$, and is generated by
all conjugates of $\beta_1^2$.

Let $\CC$ be a $k$-linear braided category over an
algebraically closed field $k$ of arbitrary characteristic.  The
braiding structure affords us representations of $\B_n$ as
follows. For any object
$X$ in $\CC$ we have braiding isomorphisms $c_{X,X}\in\End(X^{\ot 2})$ so that defining
$$
\check{R}_i:=\Id_X^{\otimes (i-1)}\otimes c_{X,X}\otimes
\Id_X^{\otimes (n-i-1)}\in\End(X^{\otimes n})
$$
we obtain a representation $\phi^n_X$ of $\B_n$
by automorphisms of $X^{\otimes n}$
by
$$\phi_X^n(\beta_i)= \check{R}_i.
$$
Similarly, for any collection of objects $\{X_i\}_{i=1}^n$,
one has representations of $\PP_n$ on $X_1\ot\cdots\ot X_n$.
Throughout the paper, when we refer to representations of $\PP_n$ and
$\B_n$ arising from tensor products of objects in a braided category,
these are the representations we mean.

We say that $Y$ is a \textbf{subobject} of $Z$ if there exists a
monomorphism $q\in\Hom_\CC(Y,Z)$, and $W$ is a \textbf{quotient
object} of $Z$ if there exists an epimorphism $p\in\Hom_\CC(Z,W)$.
Because of the functoriality of the braiding,
we have the following obvious lemma, which will be used in Section \ref{mt}.

\begin{lemma}\label{subobject}
(i) If $Y$ is a quotient object or a subobject of $Z$, then
$\phi_Y^n(\B_n)$ is a quotient group of $\phi_Z^n(\B_n)$ and similarly for the restrictions of
these representations to $\PP_n$.

(ii) Let $S$ be a finite set of objects of a braided tensor category
$\mathcal C$
for which the image of the representation of ${\mathcal{P}}_n$ in
$\End(X_1\ot \cdots\ot X_n)$ is finite for all $X_1,\ldots,X_n\in S$.
Let $X$ be the direct sum of finitely many objects taken from $S$.
Then the image of the representation of ${\mathcal {B}}_n$ in
$\End(X^{\ot n})$ is finite.
\end{lemma}

\section{The twisted quantum double of a finite group}\label{dg}
In this section we define the twisted quantum double of a finite
group, and give some basic results that we need.
For more details, see for example \cite{CP,D,witherspoon96}.

Let $k$ be an algebraically closed field of arbitrary characteristic $\ell$.
Let $G$ be a finite group with identity element $e$,
$kG$ the corresponding  group algebra, and $(k G)^*$ the dual algebra of
linear functions from $k G$ to $k$, under pointwise multiplication.
There is a basis of $(k G)^*$ consisting of the dual functions
$\delta_g$ ($g\in G$), defined by $\delta_g(h)=\delta_{g,h}$
($g,h\in G$).
Let $\omega: G\times G\times G\rightarrow k^{\times}$
be a 3-cocycle, that is
$$
  \omega(a,b,c)\omega(a,bc, d)\omega(b,c, d)=
  \omega(ab,c,d)\omega(a,b,cd)
$$
for all $a,b,c,d\in G$.
The {\bf{twisted quantum}} (or {\bf {Drinfeld}}) {\bf {double}} $D^{\omega}(G)$
is a quasi-Hopf algebra whose underlying
vector space is $(k G)^*\ot k G$.
We abbreviate the basis element $\delta_x\ot g$ of $D^{\omega}(G)$
by $\delta_x\overline{g}$ ($x,g\in G$).
Multiplication on $D^{\omega}(G)$ is defined by
$$
  (\delta_x\overline{g})(\delta_y\overline{h}) = \theta_x(g,h)
               \delta_{x,gyg^{-1}} \delta_x \overline{gh},
$$
where
$$
  \theta_x(g,h)=\frac{\omega(x,g,h)\omega(h,h,h^{-1}g^{-1}xgh)}
                {\omega(g,g^{-1}xg,h)}.
$$
As an algebra, $D^{\omega}(G)$ is semisimple if and only if the
characteristic $\ell$ of $k$ does not divide the order of $G$ \cite{witherspoon96}.

The  quasi-coassociative coproduct $\Delta: D^{\omega}(G)\rightarrow
D^{\omega}(G)\ot D^{\omega}(G)$ is defined by
$$
  \Delta(\delta_x\overline{g})=\sum_{\substack{y,z\in G\\yz=x}} \gamma_g
                (y,z) \delta_y\overline{g}\ot \delta_x\overline{g},
$$
where
$$
  \gamma_g(y,z)=\frac{\omega(y,z,g)\omega(g,g^{-1}yg,g^{-1}zg)}
     {\omega(y,g,g^{-1}zg)}.
$$
The quasi-Hopf algebra $D^{\omega}(G)$ is quasitriangular with
$$
   R=\sum_{g\in G}\delta_g\ot \overline{g} \ \ \ \mbox{ and } \ \ \
   R^{-1}=\sum_{g,h\in G}\theta_{ghg^{-1}}(g,g^{-1})^{-1} \delta_g\overline{e}
             \ot \delta_h\overline{g^{-1}}.
$$
In particular $R\Delta(a)R^{-1}=\sigma(\Delta(a))$ for all $a\in D^{\omega}(G)$,
where $\sigma$ is the transposition map.
If $X$ and $Y$ are $D^{\omega}(G)$-modules, then $\check{R}=\sigma\circ R$ provides
a $D^{\omega}(G)$-module isomorphism from $X\ot Y$ to $Y\ot X$.
Let $c_{X,Y}$ be this action by $\check{R}$.
Then the category $\Mod$-$D^{\omega}(G)$ of finite dimensional $D^{\omega}(G)$-modules
is a braided category with braiding $c$.

\section{The images of ${\mathcal {B}_n}$ and ${\mathcal{P}_n}$}\label{mt}

In this section we fix a finite group $G$ and a 3-cocycle $\omega$, and prove that the image of ${\mathcal {B}_n}$
in $\End_{D^{\omega}(G)}(V^{\ot n})$ is finite for any positive integer $n$ and
any finite dimensional $D^{\omega}(G)$-module $V$.
In case $G$ is a $p$-group, we prove that the image of ${\mathcal{P}_n}$
in $\End_{D^{\omega}(G)}(V^{\ot n})$ is also a $p$-group.

\begin{rmk} It follows from a theorem of C.\ Vafa (see \cite[Theorem 3.1.19]{BK})
and the so-called \emph{balancing axioms} that for braided fusion
categories over $\C$, the images
of the braid group generators $\beta_i$ in the above representations of $\B_n$ always have finite order.  This is far from enough to conclude that the image of $\B_n$ is finite; Coxeter \cite{Cox} has shown that the quotient of $\B_n$ by the normal closure of the subgroup generated by $\{ \beta_i^k: 1\leq i\leq n-1 \}$ is finite if and only if $\frac{1}{n}+\frac{1}{k}>\frac{1}{2}$.
\end{rmk}

\subsection*{The case of general finite groups}

Let $r$ and $m$ be positive integers. The {\textbf {full monomial group}} $G(r,1,m)$ is
the multiplicative group consisting of the $m\times m$
matrices having exactly one nonzero entry in each row and column, all of whose
nonzero entries are $r$th roots of unity.
It is one of the irreducible complex reflection groups.

Let $r=|G|_{\ell'}$ be the part
of $|G|$ not divisible by the characteristic $\ell$ of $k$
(i.e. $|G|=r\ell^s$ and $(r,\ell)=1$).

\begin{thm}\label{main}
Let $V$ be a finite dimensional $D^{\omega}(G)$-module.
Then the image of ${\mathcal{B}}_n$ in $\End(V^{\ot n})$ is
finite. More specifically, this image is a quotient of a subgroup
of $G(r,1,m)$ for $m=|G|^{2n}$.
\end{thm}

\begin{proof}
We will need the following well known lemma,
which follows from \cite[Theorem 6.5.8]{We}.
Let $\mu_r\subset k^\times$ be the set of $r$-th
roots of unity.

\begin{lemma}\label{orde}
The natural map $H^i(G,\mu_r)\to H^i(G,k^\times)$
is surjective. In particular,
any element in $H^i(G,k^\times)$
may be represented by a cocycle taking values in $\mu_r$.
\end{lemma}

Now we turn to the proof of the theorem.
As any finite dimensional $D^{\omega}(G)$-module is finitely generated,
and therefore is a quotient of a finite rank free module,
by Lemma \ref{subobject} (i),
it suffices to prove the statement when $V$ is a finite rank free module.
By Lemma \ref{subobject} (ii), we need only
consider the case $V=D^{\omega}(G)$, the left regular module.

Assume first that $n=2$.
Let $x,y,a,b\in G$.
The action of $\check{R}$ on the basis element $\delta_x\overline{a}
\ot \delta_y\overline{b}$ of $D^{\omega}(G)\ot D^{\omega}(G)$ is
\begin{eqnarray*}
  \check{R}(\delta_x\overline{a}\ot\delta_y\overline{b})&=&
   \sigma (\sum_{g\in G}\delta_g\ot \overline{g})(\delta_x\overline{a}
   \ot\delta_y\overline{b})\\
   &=& \sigma (\theta_{xyx^{-1}}(x,b) \delta_x\overline{a}\ot
         \delta_{xyx^{-1}}\overline{xb})\\
   &=& \theta_{xyx^{-1}}(x,b)\delta_{xyx^{-1}}\overline{xb}\ot\delta_x
               \overline{a}.
\end{eqnarray*}

If $n>2$, similar calculations show that each $\check{R_i}$ permutes the
chosen basis of $D^{\omega}(G)$ up to scalar multiples of the form
$\theta_{xyx^{-1}}(x,b)$.
By Lemma \ref{orde}, whe may assume that $\omega$ and hence
$\theta$ takes values in the $r$-th roots of unity.
This implies that the image of ${\mathcal{B}}_n$ in
$\End(D^{\omega}(G)^{\ot n})$ is contained in
$G(r,1,m)$.
\end{proof}

\begin{corollary} \label{grpth}
Let $\CC$ be a braided fusion category that is group-theoretical in the
sense of \cite{ENO}. Let $V$ be any object of $\CC$. Then the image of $\B_n$ in
$\End(V^{\otimes n})$ is finite.
\end{corollary}
\begin{proof} Let $Z(\CC)$ be the Drinfeld center of $\CC$. Since $\CC$ is braided, we have
a canonical braided tensor functor $F: \CC\rightarrow Z(\CC)$. Thus it suffices
to show  the result holds for the category $Z(\CC)$. Since $\CC$ is
group-theoretical, $Z(\CC)$ is equivalent to $\Mod\mbox{-}D^\omega(G)$ for some $G$, $\omega$. Thus the desired
result follows from Theorem \ref{main}.
\end{proof}

\subsection*{The case of $p$-groups.}

\begin{theorem}\label{pgroup}
Suppose that $G$ is a finite $p$-group and $V$ is a finite dimensional $D^{\omega}(G)$-module.
Then the image of ${\mathcal{P}}_n$ in $\End(V^{\ot n})$ is also a $p$-group.
\end{theorem}

The rest of the subsection is occupied by the proof of
Theorem \ref{pgroup}. We will need a technical lemma:

\begin{lemma}\label{plemma2}
Let $H$ be a group with normal subgroups
$H=H_0\supset H_1\supset ...\supset H_N=1$, such that $H_i/H_{i+1}$ is
abelian, and $[H_i,H_j]\subset H_{i+j}$, and let $I$ be a subgroup of
$\Aut(H)$ that preserves this filtration and acts trivially
on the associated graded group. Then $I$ is nilpotent of class
at most $N-1$.
\end{lemma}
\begin{proof}
 Let $L_1(I)=I$, $L_2(I)=[I,I]$, $L_3(I)=[[I,I],I],\ldots$, be the lower central series
of $I$. We must show $L_{N}(I)=1$.

We prove by induction on $n$ that for any $f\in L_n(I)$ and $h\in H_m$,
$f(h)=ha(h)$, where $a(h)\in H_{n+m}$.

The  case $n=1$ is clear, since
$f\in I$ acts trivially on $H_m/H_{m+1}$.
Suppose  the statement is true for $n$.
Take $g\in I$, $f\in L_n(I)$ and $h\in H_m$ so that:
$f(h)=ha(h)$, $g(h)=hb(h)$, where $a(h)\in H_{n+m}$ and $b(h)\in H_{m+1}$.
Then
$fg(h)=f(h)f(b(h))=ha(h)b(h)a(b(h))$, while
$gf(h)=hb(h)a(h)b(a(h))$.

Since $g$ acts trivially on the associated graded group,
$b(a(h))\in H_{n+m+1}$. Also $a(b(h))\in H_{n+m+1}$
since $b(h)\in H_{m+1}$, by the induction assumption.
Moreover, $a(h)b(h)=b(h)a(h)$ modulo $H_{n+m+1}$ since $[H_i,H_j]\subset H_{i+j}$.
Thus, $fg(h)=gf(h)$ in $H/H_{n+m+1}$, and thus $[f,g](h)=h$ in
$H/H_{n+m+1}$, which is what we needed to show.

Taking $m=0$ and $n=N-1$, any $[f,g]\in L_{N}(I)$
is the identity on $H=H/H_{N}$, and the lemma is proved.
\end{proof}

Now we are ready to prove the theorem.
Any finite dimensional $D^{\omega}(G)$-module is a quotient of a
multiple of the left regular $D^{\omega}(G)$-module $H=D^{\omega}(G)$.
By Lemma \ref{subobject}, it suffices to show that the image of $\PP_n$ in
$\End(H^{\ot n})$ is a $p$-group.  By Theorem \ref{main}, the
image $K$ of $\PP_n$ is a subgroup of the
full monomial group $G(r,1,m)$, where $r=p^t$ for some $t$, and
$m=|G|^{2n}$.
The normal subgroup of diagonal matrices in $K$ is thus a
$p$-group,
so it is enough to show that $K$ modulo the diagonal matrices is
a $p$-group. Thus it suffices to assume that $\omega=1$ and
$H=D(G)$.

Computing, we have:
$$
\check{R}(\overline{a}\delta_x\ot \overline{b}\delta_y)   =
    \sigma(\sum_{g\in G}\delta_g\ot \overline{g})(\overline{a}\delta_x\ot \overline{b}\delta_y)
   = \overline{axa^{-1}b}\delta_y\ot \overline{a}\delta_x
$$
for all $a,b,x,y\in G$.
Denote by $(g,x)$ the element $\overline{g}\delta_x$ so that a basis of $H^{\ot n}$ is:
$$(g_1,x_1)\ot\cdots\ot(g_n,x_n)$$
with $g_i,x_i\in G$.
The braid generator $\beta_i$ fixes all factors other than the
$i$th and $(i+1)$st, and on these it acts by:
\begin{eqnarray*}
  (g_i,x_i)\ot(g_{i+1},x_{i+1})&\mapsto &
    (g_ix_ig_i^{-1} x_{i+1}, x_{i+1})\ot (g_i,x_i)\\
      &=& ([g_i,x_i]x_ig_{i+1} , x_{i+1})\ot (g_i,x_i),
\end{eqnarray*}
where $[a,b]$ denotes the group commutator.  This action induces a
homomorphism $\psi:\B_n\rightarrow \Aut(\Fr_{2n})$ where $\Fr_{2n}$ is the free group on
$2n$ generators.  Explicitly, $\psi(\beta_i)$ is the automorphism defined on generators $\{g_i,x_i\}_{i=1}^n$ of $\Fr_{2n}$ by:
\begin{eqnarray*}
&&x_j\mapsto x_j,\quad \ g_j\mapsto g_j\quad \text{for}\quad j\not\in\{i,i+1\}\\
&&x_i\mapsto x_{i+1},\quad \ x_{i+1}\mapsto x_i\\
&&g_i\mapsto [g_i,x_i]x_ig_{i+1},\quad \ g_{i+1}\mapsto g_i.
\end{eqnarray*}
Since $G$ is a $p$-group, it is nilpotent of class, say, $N-1$.
Note that $\psi$ descends to a homomorphism
$\psi_N:\B_n\rightarrow \Aut(\Fr_{2n}/L_N(\Fr_{2n}))$ where $L_N(\Fr_{2n})$ denotes the
$N$th term of the lower central series of $\Fr_{2n}$.
Since $G$ is nilpotent of class $N-1$, the action of $\B_n$ on the set $G^{2n}$ defined above factors through $\psi_N$.  Thus, setting $I=\psi_N(\PP_n)$, ones sees that the action of $\PP_n$ on $H^{\ot n}$ factors through $I$, that is, $K$ is a quotient of $I$.

Let us now show that $I$ is nilpotent.
Define a descending filtration on $M=\Fr_{2n}/L_N(\Fr_{2n})$ by
positive integers as follows. Let $M_1=M$.  Define degrees on the
generators by $\deg(g_i)=1$ and $\deg(x_i)=2$ for all $i$, and
define $M_j$ for $j\geq 2$ to be the normal closure of the group
generated by $[M_k,M_{j-k}]$ for all $0\leq k\leq j$ together
with the generators of degree at least $j$.  Since $M$ is
nilpotent, this filtration is finite. Further, $I$ preserves this
filtration and acts trivially on the quotients $M_i/M_{i+1}$.  By
Lemma \ref{plemma2}, $I$ is nilpotent.

It follows that the finite group $K$ is nilpotent.
However, $K$ is generated by conjugates of $\beta_1^2$, and we
claim that $\beta_1^2$ is an element of order a power of
$p$. Indeed, this follows from the fact that if the ground field
is $\C$ (which may be assumed without loss of generality, since the double of
$G$ is defined over the integers), then the eigenvalues of
$c_{X,Y}c_{Y,X}$ for any objects $X,Y$ are ratios of twists, which are computed from
characters of $G$ (in \cite{BK}), and hence are roots of unity of degree a power of $p$.
Therefore, $K$ is a finite $p$-group. The theorem is proved.

\section{Questions}
We mention some directions for further investigation suggested by these (and other) results.
We refer the reader to \cite{ENO} and \cite{Tur1} for the relevant definitions.

\begin{enumerate}
\item Suppose $G$ is a $p$-group.
Theorem \ref{pgroup} shows that the image of the associated
representation of $\PP_n$ is also a $p$-group.  What is its
nilpotency class relative to that of $G$? Some upper bounds can be obtained from the
proof of Theorem \ref{pgroup}, but it is not clear how tight they
are.

\item The finite groups that appear as images of representations of $\B_n$ associated to
quantum groups and link invariants at roots of unity (see \cite{FRW,GJ,J86,J89,LR,LRW}) basically fall into two classes: symplectic groups and extensions
of $p$-groups by the symmetric group $S_n$.  Does this hold for
the representations of $\B_n$ associated with $\Mod\mbox{-}D^{\omega}(G)$?  In general, is there a relationship between the image of $\PP_n$ and $G$?

\item As a modular category, $\Mod\mbox{-}D^{\omega}(G)$ gives rise to (projective) representations of
mapping class groups of compact surfaces with boundary.  Are the images always finite?  It is known to be true for the mapping class groups of the torus and the $n$-punctured sphere (Theorem \ref{main}).  For more general modular categories, the answer is definitely ``no," see \cite[Conjecture 2.4]{AMU}.
\item Let us say that a braided category $\CC$ has property $\mathcal{F}$ if all braid group representations associated to $\CC$ have finite images.  What class of braided categories have property $\mathcal{F}$?  Among braided fusion categories, Corollary \ref{grpth} shows that all braided group-theoretical categories (in the sense of \cite{ENO}) have property $\mathcal{F}$.  Do all braided fusion categories with integer Frobenius-Perron dimension have property $\mathcal{F}$?
\end{enumerate}

\end{document}